# This issue is dedicated in honour of Prof. K. Ramachandra who passed away in early 2011

# On the Half Line: K. Ramachandra

Nilotpal Kanti Sinha

*Great Lakes Institute of Management, Chennai, India*

E-mail: nilotpalsinha@gmail.com; nilotpal.sinha@greatlakes.edu.in

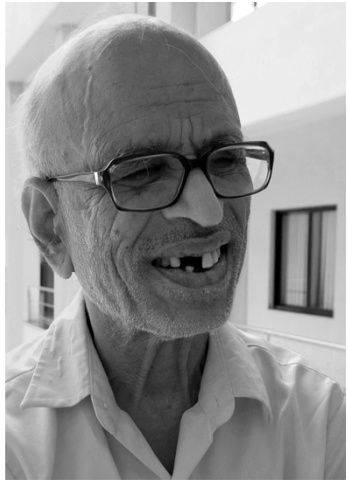

**Abstract.** A short biographical note on the life and works of K. Ramachandra, one of the leading mathematicians in the field of analytic number theory in the second half of the twentieth century.

## 1. Introduction

Kanakanahalli Ramachandra (1933–2011) was perhaps the real successor of Srinivasa Ramanujan in contemporary Indian mathematics. Ramachandra has made invaluable contributions to algebraic number theory, transcendental number theory and the theory of the Riemann zeta function. This article is a brief exposition on the life and works of Ramachandra. The title of this biographical note is motivated by the fact that Ramachandra was one of the few mathematicians who was still working on certain classical problems in number theory and many of his best results are theorems related to the values of the Riemann zeta function on the half line.

## 2. Early Life (1933–57)

Ramachandra was born on 18 August 1933 in the state of Mysore (now known as Karnataka) in southern India. His grandfather walked nearly a hundred and fifty kilometres to see the new born Ramachandra. Ramachandra hailed from a family with a modest background; his father passed away when Ramachandra was only 13. Ramachandra's mother managed his education by taking loan against their agricultural property. When Ramachandra was a student, he won a competition and was awarded a short biography of the legendary Indian mathematician Srinivasa Ramanujan. This was the book that ignited the interest for mathematics in Ramachandra.

Ramanujan's taxicab number $1729 = 9^3 + 10^3 = 1^3 + 12^3$ has become a part of mathematics folklore. During his college days, Ramachandra had a similar encounter with the number 3435. His college principal had a car with the number 3430 on the number plate. Ramachandra worked on the mathematical possibilities of this number and in the process he found that upon adding 5, the number 3435 is the only number with the unique property that when each digit was raised to a power equal to itself and the resulting numbers were added up, the sum equals the original number, i.e.

$$3^3 + 4^4 + 3^3 + 5^5 = 3435.$$

Ramachandra completed his graduation and post graduation from Central College, Bangalore. Due to family responsibilities, he had to look for a job at a young age and, just like Ramanujan, he also worked as a clerk. Ramachandra worked as a clerk at the Minerva Mills where Ramachandra's father had also worked. In spite of taking up a job quite remote from mathematics, Ramachandra studied number theory all by himself in his free time; especially the works of Ramanujan. Later, he worked as a lecturer in BMS College of Engineering. Ramachandra also served a very short stint of only six days as a teacher in the Indian Institute of science, Bangalore.

After the death of Ramanujan, Hardy wrote a series of lectures on the works of Ramanujan. These lectures were published as a book *Ramanujan: Twelve Lectures on Subjects Suggested by His Life and Work*. This became one of the classic books on Ramanujan. When Ramachandra was a college student, he found a copy of this book in the public library of the state of Mysore and he studied the works of Ramanujan with great devotion. It was this book that inspired Ramachandra to become a mathematician and in particular, a number theorist. Unfortunately someone tore the cover of the book and took away the picture of Ramanujan. This made Ramachandra very disappointed and half a century later he still did not forget the irresponsible person who tore the particular page. Later Ramachandra got a copy of the book which he always kept with him as a personal favourite.

### 3. TIFR Bombay (1958–1995)

In 1958, Ramachandra secured a post in TIFR and it was here that he met K. Chandrashekharan who was one of the experts in the theory of the Riemann zeta function in India at that time. Ramachandra studied the theory of the Riemann zeta function under K. Chandrashekharan. Later Ramachandra would himself become one of the leading experts in the theory of the Riemann zeta function and make several invaluable contributions to the subject. For the next three decades until his retirement, Ramachandra remained at TIFR and established one of the most prestigious schools of analytic number theory in collaboration with his gifted students as well as leading number theorists all over the world.

Ramachandra believed that as a mathematician one not only has to contribute to the subject but also guide the next generation of mathematicians. He worked hard to perpetuate Number theory as an active research area and succeeded in inspiring the interested students to take up the subject. Ramachandra acted as the doctoral advisor for eight students; today some of his students are among the most renowned mathematicians in the field of analytic number theory.

At the invitation of Norwegian mathematician Atle Selberg, Ramachandra went to the Institute for Advanced Study in Princeton, USA, as a visiting professor and spent a period of six months. This was Ramachandra's first foreign trip and years later when Ramachandra constructed his house in Bangalore, he named it 'Selberg House' in honor of Atle Selberg. Over the course of his career, Ramachandra visited several countries and collaborated with some of the leading number theorists of the world and also invited many of the leading mathematicians to TIFR, including the legendary Paul Erdös. Erdös visited India in 1976 and stayed as a guest in Ramachandra's house. Ramachandra is one of the few mathematicians with Erdös Number 1. He published two joint papers with Erdös.

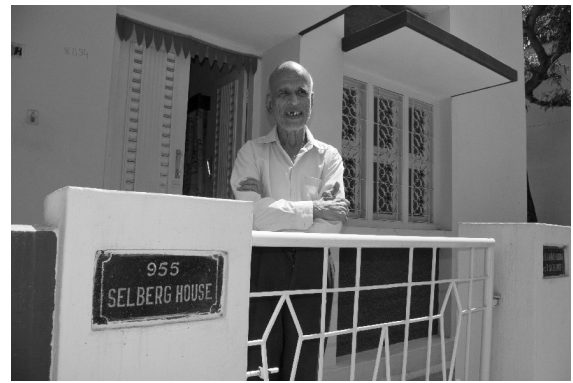

In 1978 Ramachandra founded the Hardy-Ramanujan Journal, which is considered among the most prestigious journals of number theory. It is one of the very few privately run mathematical journals in the world, funded entirely by Ramachandra and R. Balasubramanian who acted as its editors until Ramachandra passed away in 2011. The journal is published every year on 22nd December, on the occasion of the birthday of Srinivasa Ramanujan.

## 4. NIAS (1995–2011)

After retiring form TIFR, Ramachandra joined the National Institute of Advanced Studies, NIAS, Bangalore as a visiting faculty on the invitation of nuclear physicist and Founder Director of NIAS, Dr. Raja Ramanna. Ramachandra remained in NIAS and continued working on the theory of the Riemann zeta function until his death. During this period he was also associated with TIFR Bangalore as a visiting faculty. In 2003, a conference was held on the occasion of the seventieth birthday of Ramachandra in TIFR Bangalore. Several mathematicians from all over the world attended the conference and celebrated the event.

Ramachandra *left* (in the words of Erdös) on 17 January 2011. His health had broken down and he had been hospitalized for about two weeks. According to his wife Mrs. Saraswati Ramachandra: *The doctors had advised him to take complete rest and not work on mathematics. But as soon as he was discharged from the hospital, he started working on a problem that had been bothering him for the past six months.* He is survived by his wife and daughter.

## 5. Ramachandra's Mathematical Gods

Although Ramachandra was born in a brahmin family and was a devout Hindu, he also had his own perception of Gods, his mathematical Gods. He used to address the great mathematicians as Gods and his reason was simple. "*Only a God can prove such a result*," he used to say. Decorating the walls of his room at NIAS, were the poster size pictures of G. H. Hardy, Srinivasa Ramanujan and I. M. Vinogradov. 'When he spoke about these mathematicians, he used to refer them as "*They are my Gods*." Sometimes when he spoke about a particular mathematical result he used to thank his mathematical God who was an expert in that field. "*With Siegel's blessings, I was able to prove some results*," Ramachandra said, referring to the German mathematician C. L. Siegel, as he spoke about his theorems in transcendental number theory.

Once when I asked Ramachandra whom he thought was the greatest mathematician, he immediately replied, *Ivan Matveyevich*. *Vinogradov* had devised an ingenious method of evaluating exponential sums which reduced the error in the prime number theorem to

$$\pi(x) = Li(x) + O(x \exp(-c_1 (\ln x)^{3/5} (\ln \ln x)^{-1/5})).$$

In the eighty years since Vinogradov had published this result, no one was able to improve upon it. Therefore Ramachandra held Vinogradov in the highest regard. In second position he named Ramanujan and in third position he named Hardy. Ramachandra met Vinogradov twice when he was invited to attend conferences held in Russia on the occasions of the eightieth and the ninetieth birthdays of I. M. Vinogradov.

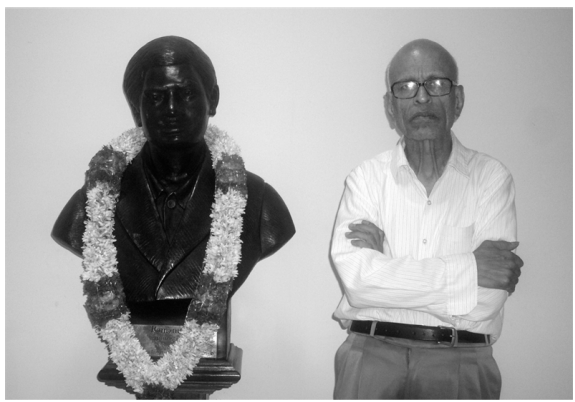

## 6. Ramachandra as a Teacher

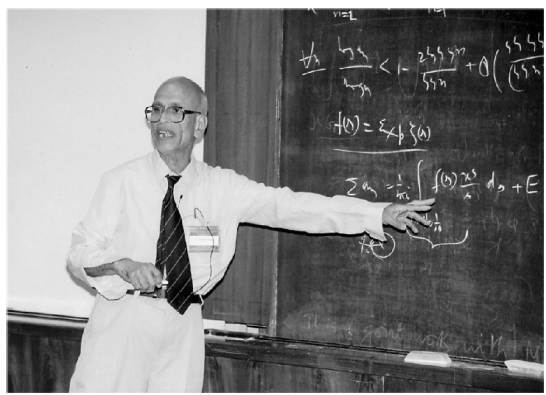

Ramachandra has mentored a few of the brightest mathematical minds today including R. Balasubramanian and T. N. Shorey. He was extremely generous to his students and he often credited his students in his papers even for very minor contributions. Ramachandra was proud of the achievements of his students and often spoke about their works. *Balu and Shorey have brought a great name to me*, he said, referring to R. Balasubramanian and T. N. Shorey. It is often said in lighter vein that Ramachandra's greatest contribution to the theory of transcendental numbers is

T. N. Shorey because Shorey has proved some extraordinary results in this field. Ramachandra often said that his student R. Balasubramanian was a better mathematician. During the days when R. Balasubramanian was his doctoral student at TIFR Bombay, Ramachandra often introduced R. Balasubramanian as: *This is R. Balasubramanian. We are teacher and student but at different times the roles will be reversed.*

### 7. Mathematical Works

During his career as a mathematician, Ramachandra published over hundred and fifty papers which included several important works in the fields of algebraic number theory, transcendental number theory and the theory of the Riemann zeta function. Ramachandra was among the pioneers in evaluating the fractional moments of the Riemann zeta function. He was also one of the first mathematicians to consider the gap between numbers with large prime factors. Several key areas of analytic number theory that Ramachandra had pioneered, continue to be active areas of research even today.

In his early years as a number theorist, Ramachandra worked in the field of algebraic number theory. His first paper was: Some applications of Kronecker's limit formula, Ann. of Math (2), 80 (1964), 104–148. The reviewer M. Eichler remarked: *This paper contains some remarkable new results on the construction of the ray class field of an imaginary quadratic number field.* Ramachandra completed his Phd under the supervision of K. G. Ramanathan at TIFR Bombay (now known as Mumbai) in 1965.

When the seminal work of Alan Baker appeared in the 1960's, Ramachandra and his students, especially T. N. Shorey, took up transcendental number theory and made remarkable contributions to both the theory and its applications to problems of classical number theory. A detailed exposition of Ramachandra's contribution can be found in Michel Waldschmidt's paper 'On Ramachandra's Contributions to Transcendental Number Theory.'

In 1974, Ramachandra turned his attention to hard core classical analytic number theory, especially the theory of Riemann zeta function and general Dirichlet series. His contributions to the theory of the Riemann zeta function is best summarized in the words of British mathematician and a Fellow of the Royal Society, Roger Heath-Brown:

*As soon as I entered research, 30 years ago, yours became a familiar name; and your infuence has remained with me ever since. Time permits me to mention in detail only one strand of your work – but it is one that clearly demonstrates how important your research has been. A little over 20 years back you proved the first results on fractional moments of the Riemann Zeta-function. At first I could not believe they were correct!! Since then however the ideas have been extended in a number of ways. They have lead of course to a range of important new results about the Zeta-function and other Dirichlet series. But just as significantly the ideas have led to new conjectures on the moments of the Riemann Zeta-function. These conjectures provide the first successful test for the application of Random Matrix Theory in this area. Nowadays this is a growing area which has contributed much to our understanding of zeta-functions. And it can all be traced back to your work in the late 1970's.*

### 8. Selected Theorems – The Little Flowers

Ramachandra used to dedicate his best results to his mathematical Gods in papers whose title began with *Little flowers to ....* For instance when he visited Russia on the occasion of the ninetieth birthday of I. M. Vinogradov, the paper that he presented in the conference was titled *Little flowers to I. M. Vinogradov*. We shall present a few flowers from Ramachandra's garden that roughly cover his genre of work. For more details on Ramachandra's work, I would request the reader to refer to the volumes of the Hardy-Ramanujan journal which are available online at http://www.imsc.res.in

**Theorem 8.1 (K. Ramachandra).** *Let $\lambda$ be any constant satisfying $1/2 < \lambda < 1$ and $l$ a non-negative integer constant. Put $H = T^\lambda$. Then we have*

$$(\ln T)^{1/4+l} \ll \frac{1}{H} \int_T^{T+H} \left| \frac{d^l}{dt^l} \zeta(1/2 + it) dt \right| \ll (\ln T)^{1/4+l}.$$

**Theorem 8.2 (R. Balasubramanian and K. Ramachandra).** *Let $t$ be a fixed transcendental number and $x \geq 1, y \geq 1$ be integers. Let $n$ be any integer such that $x \leq n < x + y$ for which $2^{t^n}$ defined as $\exp(t^n \ln 2)$ is algebraic. The number of such integers is less than or equal to $(2y)^{1/2} + O(y^{1/4})$.*

**Theorem 8.3 (M. Jutila, K. Ramachandra and T. N. Shorey).** *Let $k > 2$ and $n_1 = n_1(k), n_2 = n_2(k), \ldots$ be the sequence of all positive integers which have at least one prime factor $> k$. Put $f(k) = \max(n_{i+1} - n_i)$ the maximum being taken over all $i > 1$. Then*

$$f(k) \ll \frac{k}{\ln k}\left(\frac{\ln \ln \ln k}{\ln \ln k}\right).$$

**Theorem 8.4 (K. Ramachandra).** *For all sufficiently large $m$, between $m^2$ and $(m+1)^2$, there is an $n$ and a prime $p$ dividing $n$ such that $p > n^{1/2+1/11}$.*

**Theorem 8.5 (Ramachandra, Shorey and Tijdeman).** *There exists an absolute constant $c_2 > 0$ such that for $n \geq 3$ and $g = g(n) = \left[c_2\left(\frac{\ln n}{\ln \ln n}\right)^3\right]$, it is possible to choose distinct primes $p_1, p_2, \ldots, p_g$ such that $p_i | (n+1)$ for $1 \leq i \leq g$.*

## 9. Acknowledgement

I am indebted to Mrs. Saraswati Ramachandra, wife of K. Ramachandra, for sharing the details of the last few days of Ramachandra's life. She was kind enough to give me a small photograph of Ramachandra that she had preserved. I would also like to thank C. S. Aravinda, TIFR Bangalore, and Kishor Bhat, NIAS Bangalore, for providing me with valuable information and encouraging me to write this biography. I wish to thank S. L. Ravi Shankar for providing some of the photographs in this article. I would also like to thank Ankita Jain and Roshani Nair, GLIM Chennai, for reviewing this article and giving valuable suggestions.

# Prof. K. Ramachandra: Reminiscences of his Friends

### M. Pavaman Murthy

*A Few Reminiscences of K. Ramachandra in his Early Years at TIFR*

I joined TIFR in 1958 along with Ramachandra, R. R. Simha and Vasanthi Rao. Ramachandra stood out in this group of four because of his remarkable dedication to mathematics and his knowledge of advanced number theory. He had made up his mind even before joining TIFR to work in number theory. A few months after joining the institute, Ramachandra was studying research papers in number theory when the rest of us in the group were struggling with topics exposed in the "baby seminars." The four of us were good friends. In recent times, whenever I visited TIFR from Chicago, Ramachandra used to enter in to my office and inquire about my family and tell me about his work and his students' work. He always spoke of his students with pride.

Here are some of my reminiscences of Ramachandra in his early years at TIFR.

In the first or second year after joining the institute, one early morning in TIFR hostel in the Old Yacht Club, I still remember witnessing a heated discussion between Ramachandra and Raghavan Narasimhan comparing the contributions of Hilbert and Kummer to number theory. Probably that was the first time that I came to know that there was a famous number theorist called Kummer.

Here is a glimpse of Ramachandra's sense of humor. In our early TIFR years, the institute was located in the Old Yacht Club building adjacent to Gateway of India. Several of the members of the School of Mathematics along with us new comers worked in a large hall with many tables. One day Ramachandra came to my table and showed me a theorem of Siegel and said "See Siegel proves this theorem in ten pages. I have proved the same in five pages." Indeed Siegel's proof covered ten pages. What Ramachandra had done was to copy verbatim Siegel's proof in his note book in five pages (with very small handwriting)!

Sometime after joining the institute, we four new comers were told that there would not be any oral examinations at the